\documentclass[12pt]{article}
\usepackage{amsmath,amssymb,amsthm}
\usepackage{fullpage}

\newtheorem{thm}{Theorem}
\newtheorem{cor}{Corollary}
\newcommand{\dis}{\displaystyle}

\begin{document}

\title{Ramanujan's Approximation to the $n\,$th Partial Sum\\
of the Harmonic Series}
\author{Mark B. Villarino\\
Depto.\ de Matem\'atica, Universidad de Costa Rica,\\
2060 San Jos\'e, Costa Rica}
\date{June 5, 2004}

\maketitle

 \begin{abstract}
 A simple integration by parts and telescopic cancellation leads to
 a derivation of the first two terms of Ramanujan's asymptotic
series for the $n$th partial sum of the harmonic series. Kummer's
 transformation gives three more terms with an explicit error
 estimate. We also give best-possible estimates of Lodge's
 approximations.
 \end{abstract}

Entry~9 of Chapter~38 of B. \textsc{Berndt}'s edition of
\textsc{Ramanujan}'s Notebooks, Volume 5 \cite[p.~521]{Berndt} reads
(in part):
\begin{quote}
``\textit{Let $m :=  \frac{n(n+1)}{2}$, where $n$ is a positive
integer. Then, as $n$ approaches infinity,}
$$
\sum_{k=1}^n \frac{1}{k}
\sim \frac{1}{2} \ln(2m) + \gamma + \frac{1}{12m} - \frac{1}{120m^2}
+ \frac{1}{630m^3} - \frac{1}{1680m^4} + \frac{1}{2310m^5} - [\cdots]
\text{.''}
$$
\end{quote}
(The entry includes terms up to $m^{-9}$.)

Berndt's proof simply verifies (as he himself explicitly notes) that
Ramanujan's expansion coincides with the standard \textsc{Euler}
expansion
\begin{align*}
H_n := \sum_{k=1}^n \frac{1}{k}
&\sim \ln n + \gamma + \frac{1}{2n} - \frac{1}{12n^2}
+ \frac{1}{120n^4} -[ \cdots]
\\
&= \ln n + \gamma - \sum_{k=1}^\infty \frac{B_k}{n^k} 
\end{align*}
where $B_k$ denotes the $k^{\mathrm{th}}$ \textsc{Bernoulli} number
and $\gamma := 0.57721\cdots$ is Euler's constant.

However, Berndt does \textit{not} show that Ramanujan's expansion is
\textit{asymptotic} in the sense that the error in the value obtained
by stopping at any particular stage in Ramanujan's series is less than
the next term in the series. Indeed we have been unable to find
\emph{any} error analysis of Ramanujan's series.

We therefore offer the following error analysis which shows \emph{the
first five terms} of Ramanujan's series to be asymptotic in the sense
above. Our methods can be extended to any number of terms in the
expansion, but we will limit our presentation to the first five.

Berndt also states that there is no ``natural'' way to obtain an
expansion of $H_n$ in powers of $m$. In fact, our method produces the
first two terms of Ramanujan's expansion, namely,
$ \frac{1}{12m} - \frac{1}{120m^2}$, \textit{automatically}, and
the later terms by a simple ``\textsc{Kummer}'s transformation'': see
\cite[p.~260]{Knopp}.

\begin{thm}
Let $m := \dfrac{1}{2}n(n+1)$, where $n$ is a positive integer. Then
there exists a $\Theta_n$, with $0 < \Theta_n < 1$, for which the
following equation is true:
$$
\fbox{$\dis
\sum_{k=1}^n \frac{1}{k}
= \frac{1}{2} \ln(2m) + \gamma + \frac{1}{12m} - \frac{1}{120m^2}
+ \frac{1}{630m^3} - \frac{1}{1680m^4} + \frac{\Theta_n}{2310m^5}.
$}
$$
\end{thm}

We observe that this expansion of $H_n$ \textit{does have} the
property that the error in the value obtained by stopping at any
particular stage in it is less than the next term since the terms
alternate in sign and decrease monotonically in absolute value.

\begin{proof}
We follow a hint from \textsc{Bromwich}
\cite[p.~460, Exercise~18]{Bromwich}: set
$$
\epsilon_n := H_n - \frac{1}{2} \ln[n(n+1)] - \gamma .
$$
Then (this is Bromwich's hint),
$$ 
\epsilon_{n-1} - \epsilon_n = \int_0^1 \frac{t^2}{n(n^2 - t^2)} \,dt .
$$
Therefore,
\begin{align*}
\epsilon_n
&= (\epsilon_n - \epsilon_{n+1}) + (\epsilon_{n+1} - \epsilon_{n+2})
+ \cdots
\\
&= \sum_{k=n+1}^\infty \int_0^1 \frac{t^2}{k(k^2 - t^2)} \,dt .
\end{align*} 
Now, integrating by parts and using the partial fraction expansion  
of $\dis\frac{1}{k(k^2 - 1)}$, and then integrating by parts again 
and using the partial fraction expansion of
$\dis\frac{1}{k(k^2 - 1)^2}$, we obtain


\noindent
\begin{align*}
\epsilon_n
&= \sum_{k=n+1}^{\infty} \biggl\{ \frac{1}{3k(k^2 - 1)}  
- \frac{2}{3} \int_0^1 \frac{t^4}{k(k^2 - t^2)^2} \,dt  \biggr\}
\\
&= \sum_{k=n+1}^{\infty} \biggl\{ \frac{1}{6}
\biggl( \frac{1}{[k-1]k} - \frac{1}{k[k+1]} \biggr)  
- \frac{2}{3} \int_0^1 \frac{t^4}{k(k^2 - t^2)^2} \,dt  \biggr\}
\\
&= \frac{1}{6n(n+1)} - \sum_{k=n+1}^{\infty} \biggl\{
\frac{2}{15k(k^2 - 1)^2} - \frac{8}{15} \int_0^1
\frac{t^6}{k(k^2 - t^2)^3} \,dt \biggr\}
\\
&\hspace{9em}
- \frac{8}{15} \int_0^1 \frac{t^6}{k(k^2 - t^2)^3} \,dt  \biggr\}
\\
&= \frac{1}{6n(n+1)} - \sum_{k=n+1}^{\infty} \biggl\{ \frac{2}{15}  
\biggl[ \frac{1}{4(k-1)^2} - \frac{1}{2k(k-1)} - \frac{1}{4(k+1)^2}
        + \frac{1}{2k(k+1)} \biggr]
\\
&\hspace{9em}
- \frac{8}{15} \int_0^1 \frac{t^6}{k(k^2 - t^2)^3} \,dt  \biggr\}
\\
&= \frac{1}{6n(n+1)} - \frac{2}{15} \biggl( -\frac{1}{2n[n+1]}
+ \frac{1}{4n^2} + \frac{1}{4[n+1]^2} \biggr)
+ \sum_{k=n+1}^{\infty} \frac{8}{15} \int_0^1
\frac{t^6}{k(k^2 - t^2)^3} \,dt      
\\
&= \frac{1}{6n(n+1)}
- \frac{2}{15} \biggl( \frac{1}{4n^2[n+1]^2} \biggr)
+ \sum_{k=n+1}^{\infty} \frac{8}{15} \int_0^1
\frac{t^6}{k(k^2 - t^2)^3} \,dt      
\\
&= \frac{1}{12m} - \frac{1}{120m^2} + \frac{8}{15}
\sum_{k=n+1}^{\infty} \int_0^1 \frac{t^6}{k(k^2 - t^2)^3} \,dt ,
\end{align*}
and we have obtained \emph{the first two terms} of Ramanujan's
expansion in powers of $m$ in a very simple and straightforward
manner.

A third integration by parts gives us
\begin{align*}
\frac{8}{15} \sum_{k=n+1}^{\infty} \int_0^1
\frac{t^6}{k(k^2 - t^2)^3} \,dt
&= \frac{8}{105} \sum_{k=n+1}^{\infty} \frac{1}{k(k^2 - 1)^3}
- \frac{16}{35} \sum_{k=n+1}^{\infty} \int_0^1
\frac{t^8}{k(k^2 - t^2)^4} \,dt .
\end{align*}

Unfortunately, the series
$$
\dis  \frac{8}{105} \sum_{k=n+1}^{\infty} \frac{1}{k(k^2 - 1)^3}
$$
apparently does not lead to a nice partial fractions telescopic
cancellation, and so we need a \textit{new idea}. If we look at the
\textit{asymptotic} behavior, as $n \to \infty$, of the error term
$$
\dis \frac{8}{15} \sum_{k=n+1}^{\infty} \int_0^1
\frac{t^6}{k(k^2 - t^2)^3} \,dt,
$$
we observe that the definition of $m$ implies that $n^2 \sim 2m$ and 
therefore the error term is
\begin{align*}
\frac{8}{15} \sum_{k=n+1}^{\infty} \int_0^1
\frac{t^6}{k(k^2 - t^2)^3} \,dt
&= \frac{8}{105} \sum_{k=n+1}^{\infty} \frac{1}{k(k^2 - 1)^3}
- \frac{16}{35} \sum_{k=n+1}^{\infty} \int_0^1
\frac{t^8}{k(k^2 - t^2)^4} \,dt
\\
&\sim \frac{8}{105} \sum_{k=n+1}^{\infty} \frac{1}{k^7}
\\
&\sim \frac{8}{105} \cdot \frac{1}{6n^6}
\\
&\sim \frac{8}{105} \cdot \frac{1}{48m^3} 
= \frac{1}{630m^3} 
\end{align*}
and the next term in Ramanujan's expansion is indeed
$\dis \frac{1}{630m^3}$.

\vspace*{6pt}

The new idea is this. We observe that the asymptotic error term
$$
\frac{1}{630m^3} = \frac{8}{105} \cdot \frac{1}{6n^3(n+1)^3}
$$
can be represented as \textit{the sum of a telescopic series} as
follows:
\begin{align*}
\frac{1}{6n^3(n+1)^3}
&= \frac{1}{6n^3(n+1)^3} - \frac{1}{6(n+1)^3(n+2)^3} + \\
&  \frac{1}{6(n+1)^3(n+2)^3} - \frac{1}{6(n+2)^3(n+3)^3} + [\cdots].                         
\end{align*}

Therefore, if we \textit{add and subtract} this expansion from the
error term series we obtain
\begin{align*}
\dis \frac{8}{105}\dis  \sum_{k=n+1}^{\infty} \frac{1}{k(k^2 - 1)^3}
&= \frac{8}{105} \frac{1}{6n^3(n+1)^3} 
+ \frac{8}{105} \sum_{k=n+1}^{\infty} \left\{ \frac{1}{k(k^2 - 1)^3}
- \left[ \frac{1}{6(k-1)^3k^3} - \frac{1}{6k^3(k+1)^3} \right] \right\}
\\
&= \frac{1}{630m^3} - \frac{8}{315} \sum_{k=n+1}^{\infty}
\frac{1}{k^3(k^2-1)^3},
\end{align*}
and the error in Ramanujan's expansion takes the form
\begin{align*}
\epsilon_n
&:= H_n-\dis\frac{1}{2}\ln (2m)-\gamma
\\
&= \frac{1}{12m}-\frac{1}{120m^2}+\frac{1}{630m^3}
- \frac{8}{315} \sum_{k=n+1}^{\infty}\frac{1}{k^3(k^2-1)^3}
- \frac{16}{35} \sum_{k=n+1}^{\infty} \int_0^1
\frac{t^8}{k(k^2 - t^2)^4} \,dt ,
\end{align*}
and we now have \textit{three} terms of Ramanujan's expansion. This
technique is a simple example of Kummer's transformation.

To extend the expansion to terms of order $m^5$, we integrate by parts
three times and then apply this Kummer transformation technique twice.
We obtain
\begin{align*}
\epsilon_n
&= \frac{1}{12m} - \frac{1}{120m^2} + \frac{1}{630m^3}
- \frac{8}{315} \sum_{k=n+1}^{\infty} \frac{1}{k^3(k^2-1)^3}
- \frac{16}{315} \sum_{k=n+1}^{\infty} \frac{1}{k(k^2-1)^4}
\\
&\qquad + \frac{128}{3465} \sum_{k=n+1}^{\infty} \frac{1}{k(k^2-1)^5}
- \frac{256}{9009} \sum_{k=n+1}^{\infty} \frac{1}{k(k^2-1)^6}
+ \frac{1024}{3003} \sum_{k=n+1}^{\infty} \int_0^1
\frac{t^{14}}{k(k^2 - t^2)^7} \,dt
\\
&= \frac{1}{12m} - \frac{1}{120m^2} + \frac{1}{630m^3}
- \frac{1}{1680m^4}
+ \frac{32}{315} \sum_{k=n+1}^{\infty} \frac{1}{k^3(k^2-1)^4}
+ \frac{128}{3465} \sum_{k=n+1}^{\infty} \frac{1}{k(k^2-1)^5}
\\
&\qquad - \frac{256}{9009} \sum_{k=n+1}^{\infty} \frac{1}{k(k^2-1)^6}
+ \frac{1024}{3003} \sum_{k=n+1}^{\infty} \int_0^1
\frac{t^{14}}{k(k^2 - t^2)^7}\,dt
\\
&= \frac{1}{12m} - \frac{1}{120m^2} + \frac{1}{630m^3}
- \frac{1}{1680m^4} + \frac{1}{2310m^5}
\\
&\qquad -\frac{32}{3465} \sum_{k=n+1}^{\infty}
\frac{41k^2+3}{k^5(k^2-1)^5}
- \frac{256}{9009} \sum_{k=n+1}^{\infty} \frac{1}{k(k^2-1)^6}
+ \frac{1024}{3003} \sum_{k=n+1}^{\infty} \int_0^1
\frac{t^{14}}{k(k^2 - t^2)^7}\,dt .
\end{align*}
To bound $\epsilon_n$ from \emph{above}, we observe that
\begin{align*}
\epsilon_n
&= \frac{1}{12m} - \frac{1}{120m^2} + \frac{1}{630m^3}
- \frac{1}{1680m^4} + \frac{1}{2310m^5}
\\
&\qquad - \Bigg\{ \frac{32}{3465} \sum_{k=n+1}^{\infty}
\frac{41k^2+3}{k^5(k^2-1)^5}
+ \frac{256}{9009} \sum_{k=n+1}^{\infty} \frac{1}{k(k^2-1)^6}
- \frac{1024}{3003} \sum_{k=n+1}^{\infty} \int_0^1
\frac{t^{14}}{k(k^2 - t^2)^7}\,dt \Bigg\} .
\end{align*}
We will show that the term in curly brackets is \emph{positive}. In
fact,
\begin{align*}
& \frac{32}{3465} \sum_{k=n+1}^{\infty} \frac{41k^2+3}{k^5(k^2-1)^5}
+ \frac{256}{9009} \sum_{k=n+1}^{\infty} \frac{1}{k(k^2-1)^6}
- \frac{1024}{3003} \sum_{k=n+1}^{\infty} \int_0^1
\frac{t^{14}}{k(k^2 - t^2)^7} \,dt
\\
&> \frac{32}{3465} \sum_{k=n+1}^{\infty}
\frac{41k^2}{k^5(k^2-1)^5}
+ \frac{256}{9009} \sum_{k=n+1}^{\infty} \frac{1}{k(k^2-1)^6}
- \frac{1024}{3003} \sum_{k=n+1}^{\infty} \frac{1}{(k-1)^{15}}
  \int_0^1 t^{14} \,dt
\\
&= \frac{1312}{3465} \sum_{k=n+1}^{\infty}
\frac{k^2}{k^5(k^2-1)^5}
+ \frac{256}{9009} \sum_{k=n+1}^{\infty} \frac{1}{k(k^2-1)^6}
- \frac{1024}{3003} \sum_{k=n+1}^{\infty} \frac{1}{15(k-1)^{15}}
\\
&> \biggl( \frac{1312}{3465} + \frac{256}{9009} \biggr)
\sum_{k=n+1}^{\infty} \frac{1}{k^{13}}
- \frac{1024}{45045} \sum_{k=n+1}^{\infty} \frac{1}{(k-1)^{15}}
\\
&= \sum_{k=n+1}^{\infty}\left[\frac{6112}{15015}\, \frac{1}{k^{13}}
- \frac{1024}{45045}\, \frac{1}{(k-1)^{15}}\right] ,
\\
\end{align*}
and a simple exercise in inequalities shows that each summand of the
rightmost sum is positive whenever $k > n \geqslant 3$. Therefore,
$$
\epsilon_n < \frac{1}{12m} - \frac{1}{120m^2} + \frac{1}{630m^3}
- \frac{1}{1680m^4} + \frac{1}{2310m^5} .
$$

To bound $\epsilon_n$ from \emph{below}, we note that
\begin{align*}
\epsilon_n
&> \frac{1}{12m} - \frac{1}{120m^2} + \frac{1}{630m^3}
- \frac{1}{1680m^4} + \frac{1}{2310m^5}
\\
&\qquad - \biggl( \frac{32}{3465} \sum_{k=n+1}^{\infty}
\frac{41k^2 + 3}{k^5(k^2-1)^5}
+ \frac{256}{9009} \sum_{k=n+1}^{\infty} \frac{1}{k(k^2-1)^6} \biggr)
\\
&> \frac{1}{12m} - \frac{1}{120m^2} + \frac{1}{630m^3}
- \frac{1}{1680m^4} + \frac{1}{2310m^5}
- \biggl( \frac{32}{3465} \cdot 42 + \frac{256}{9009} \biggr)
 \sum_{k=n+1}^{\infty} \frac{1}{(k-1)^{13}}
\\
&> \frac{1}{12m} - \frac{1}{120m^2} + \frac{1}{630m^3}
- \frac{1}{1680m^4} + \biggl\{ \frac{1}{2310m^5} 
- \frac{4688}{135135(n-\frac{1}{2})^{12}} \biggr\}
\\
&> \frac{1}{12m} - \frac{1}{120m^2} + \frac{1}{630m^3}
- \frac{1}{1680m^4},
\end{align*}
since, as another simple exercise in inequalities shows, the term in
curly brackets is positive for $n \geqslant 5$. Finally, the theorem
can be checked directly for $n = 1,2,3,4$.
\end{proof}

Now we note two fascinating corolaries due, in concept, but without
the optimal error estimates, to the British mathematician
\textsc{Alfred Lodge}~\cite{Lodge}.

\begin{cor}
For every positive integer $n$, there is a $\lambda_n$ for which
$$
\fbox{$\dis 1 + \frac{1}{2} + \frac{1}{3} +\cdots+ \frac{1}{n}
= \frac{1}{2}\ln(2m) + \gamma + \frac{1}{12m+\frac{6}{5}} + \lambda_n$}
$$
where
$$
0 < \lambda_n < \frac{19}{25200m^3}.
$$  
In fact,
$$
\lambda_n = \frac{19}{25200m^3} - \rho_n ,
$$ 
where $0 < \rho_n < \dis\frac{43}{84000m^4}$. The constants
$\dis \frac{19}{25200}$ and $\dis \frac{43}{84000}$ are the best
possible.
\qed
\end{cor}

\vspace{6pt}

\begin{cor}
For every positive integer $n$, define the quantity $\Lambda_n$ by
the following equation:
$$
1 + \frac{1}{2} + \frac{1}{3} +\cdots+ \frac{1}{n}
=: \frac{1}{2}\ln(2m) + \gamma + \frac{1}{12m+\Lambda_n}.
$$                                                                                                                                                                       
Then
$$
\Lambda_n = \dis\frac{6}{5} - \frac{19}{175m} + \frac{13}{250m^2}
 - \frac{\delta_n}{m^3} , 
$$
where $\dis 0 < \delta_n < \dis\frac{187969}{4042500}$. The
constants in the expansion of $\Lambda_n$ all are the best possible.
\qed
\end{cor}

\vspace{6pt}

\begin{cor}
For every positive integer $n\geqslant 1$ there exists a number $c_n$,
$0 < c_n < 1$, such that the following approximation is valid:
$$
H_n = \dis \frac{1}{2}\ln(2m) + \gamma + \frac{c_n}{12m}.
\eqno \qed
$$
\end{cor}

The first and second corollaries appeared, in much less precise form
and \emph{with no error estimates}, in a very interesting paper by
Lodge \cite {Lodge}, which later mathematicians inexplicably ignored.
Lodge gives some numerical examples of the error in the approximative
equation
$$H_n \approx \frac{1}{2}\ln(2m) + \gamma + \frac{1}{12m + \frac{6}{5}}$$
in Corollary 1; he also presents the first two terms of
$\Lambda_n$ from Corollary 2. An asymptotic error estimate
for Corollary 1 (with the incorrect constant $\frac{1}{150}$
instead of $\frac{1}{165\frac{15}{19}}$) appears as Exercise 19 on
page 460 in Bromwich~\cite{Bromwich}.

Our third corollary is the exercise (no.~18, page 460) that
Bromwich originally proposed and is, of course, a trivial
consequence of our main theorem. In fact, it is due to
\textsc{E. Ces\`aro} \cite{Cesaro} who proved it in~1885 by a
completely different technique. By the way, this was two years before
Ramanujan was born!

\subsubsection*{Acknowledgment}
I thank Joseph C. V\'arilly for comments on an earlier version. 
Support from the Vicerrector\'{\i}a de Investigaci\'on of the 
University of Costa Rica is acknowledged.

\end{document}